\documentclass[twoside,a4paper,10pt]{article}
\usepackage{latexsym,amssymb}
\usepackage{mltex}
\usepackage{graphicx}
\usepackage[frenchb]{babel}
\usepackage[matrix,arrow]{xy}
\def \be{\begin{eqnarray*}}
\def \ee{\end{eqnarray*}}
\def \ben{\begin{enumerate}}
\def \een{\end{enumerate}}
\def \beit{\begin{itemize}}
\def \eeit{\end{itemize}}
\def \bui#1#2{\mathrel{\mathop{\kern 0pt#1}\limits^{#2}}}
\def \buil#1#2{\mathrel{\mathop{\kern 0pt#1}\limits_{#2}}}
\def \bfll{\begin{flushleft}}
\def \efll{\end{flushleft}}
\def \bflr{\begin{flushright}}
\def \eflr{\end{flushright}}

\def \findemo{\bflr$\Box$\eflr}
\def \lra{\longrightarrow}

\def \wih{\widehat}
\def \wit{\widetilde}
\def \wnabla{\wit{\nabla}}
\def \cdotM{\buil{\cdot}{M}}

\newcommand{\pa}[1]{\left(#1\right)}

\newtheorem{thm}{Th\'eor\`eme}

\title{Une nouvelle estimation extrins\`eque du spectre de l'op\'erateur de Dirac}
\author{Nicolas Ginoux}

\begin{document}
\maketitle

\noindent{\bf R\'esum\'e:} Nous \'etablissons une nouvelle majoration optimale pour les plus petites valeurs propres de l'op\'erateur de Dirac sur une hypersurface compacte de l'espace hyperbolique.\\

\begin{center}{\Large A new extrinsic estimate for the spectrum of the Dirac operator} \end{center}
$ $\\
\noindent{\bf Abstract:} We prove a new upper bound for the smallest eigenvalues of the Dirac operator on a compact hypersurface of the hyperbolic space.\\
$ $\\

Soit $(M^m,g)$ une hypersurface riemannienne compacte orient\'ee (de dimension $m$) d'une vari\'et\'e riemannienne spinorielle $(\wit{M}^{m+1},g)$. Notons $\lambda_1$ la plus petite valeur propre de l'op\'erateur de Dirac associ\'e \`a $g$ et \`a la structure spinorielle induite sur $M$, et $H$ la courbure moyenne de $M$ dans $\wit{M}$. Dans \cite{Baer98}, C. B\"ar montre que
\begin{eqnarray}
\lambda_1^2&\leq&\frac{m^2}{4\mathrm{Vol}(M)}\int_MH^2 v_g\,,\qquad\quad\;\textrm{ si }\wit{M}=\mathbb{R}^{m+1},\\
\lambda_1^2&\leq&\frac{m^2}{4\mathrm{Vol}(M)}\int_M(H^2+1) v_g\,,\quad\textrm{ si }\wit{M}=S^{m+1},\qquad\textrm{ et}\\
|\lambda_1|&\leq&\frac{m}{2}(\buil{\sup}{M}\hspace{-1mm}|H|+1)\,,\qquad\qquad\,\,\,\,\textrm{ si }\wit{M}=\mathbb{H}^{m+1}.\end{eqnarray}
 L'auteur montre \'egalement que, si $M$ est une sph\`ere g\'eod\'esique, alors (1) et (2) sont des \'egalit\'es, mais pas (3). Ces r\'esultats appellent donc la question suivante: peut-on am\'eliorer (3) en une majoration optimale?\\
 
\noindent Nous montrons dans cette note une estimation optimale de $\lambda_1$ pour une hypersurface de l'espace hyperbolique. La preuve de ce r\'esultat utilise une approche dif\mbox{}f\'erente de celle de \cite{Gin2002}, qui est bas\'ee sur les propri\'et\'es \emph{conformes} des op\'erateurs de Dirac et de Penrose. Nous en discutons ensuite les limites, en remarquant qu'en dimension $2$ on ne peut esp\'erer obtenir pour l'op\'erateur de Dirac une majoration $L^2$ analogue \`a celle de \cite{ElSI92} (Th\'eor\`eme 2).\\

\noindent Ce travail, qui repose sur la th\`ese de l'auteur (\cite{Ginthese}, Chap. 2), a \'et\'e ef\mbox{}fectu\'e \`a l'Institut Max-Planck pour les Math\'ematiques dans les Sciences de Leipzig, que l'auteur tient \`a remercier pour son soutien et son hospitalit\'e.

\section{Pr\'eliminaires}

\noindent Pour les pr\'eliminaires sur la g\'eom\'etrie spinorielle, on se reportera par exemple \`a \cite{LM,BHMM,BFGK,Friedlivre}.\\

 \noindent Nous rappelons quelques faits de base sur la restriction de spineurs \`a une hypersurface (on pourra aussi consulter \cite{Baer98,Morel2001,Gin2002}). Consid\'erons une hypersurface riemannienne compacte (connexe) et orient\'ee $M^m$ dans une vari\'et\'e riemannienne spinorielle $(\wit{M}^{m+1},g)$. L'orientation de $M$ permet, gr\^ace \`a l'existence sur $M$ d'un champ normal unitaire $\nu$ compatible avec les orientations de $M$ et de $\wit{M}$, d'induire une structure spinorielle sur $M$ \`a partir de celle de $\wit{M}$, poss\'edant les propri\'et\'es suivantes: notons $\Sigma M$ (resp. $\Sigma\wit{M}$) le f\mbox{}ibr\'e des spineurs de $(M,g)$ (resp. de $(\wit{M},g)$), ``$<\cdot\,,\,\cdot>$'' le produit scalaire hermitien naturel de $\Sigma M$, $\nabla$ (resp. $\wnabla$) la d\'eriv\'ee covariante canonique sur $\Sigma M$ (resp. sur $\Sigma\wit{M}$), et ``$\cdotM$'' (resp. ``$\,\cdot\,$'') la multiplication de Clif\mbox{}ford de $TM$ sur $\Sigma M$ (resp. de $T\wit{M}$ sur $\Sigma\wit{M}$).

\noindent Il existe alors un isomorphisme
\begin{equation}\label{isomspin}\Sigma\wit{M}_{|_M}\lra\left\{\begin{array}{ll}\Sigma M& \textrm{si }m\textrm{ est pair}\\\Sigma M\oplus\Sigma M& \textrm{si }m\textrm{ est impair}\end{array}\right.\end{equation} qui est unitaire et qui envoie, pour tout champ $X$ tangent \`a $M$ et toute section $\phi$ de $\Sigma\wit{M}_{|_M}$, la section $X\cdot\nu\cdot\phi$ sur $X\hspace{-1mm}\cdotM\hspace{-0.5mm}\phi$ et la section $\wnabla_X\phi$ sur $\nabla_X\phi+\frac{A(X)}{2}\hspace{-1mm}\cdotM\hspace{-0.5mm}\phi$, o\`u $A$ est le champ d'endomorphismes de Weingarten de $TM$. L'isomorphisme (\ref{isomspin}) \'etant d\'esormais assimil\'e \`a l'application identit\'e, consid\'erons l'op\'erateur $D$ agissant sur les sections de $\Sigma\wit{M}_{|_M}$ par \[D:=\left\{\begin{array}{ll}D_M& \textrm{si }m\textrm{ est pair}\\D_M\oplus -D_M& \textrm{si }m\textrm{ est impair,}\end{array}\right.\] o\`u $D_M$ est l'op\'erateur de Dirac (dit \emph{fondamental}) de $(M,g)$. Par d\'ef\mbox{}inition, $D$ est elliptique autoadjoint et $\textrm{Spec}(D^2)=\textrm{Spec}(D_M^2)$. Nous noterons, en tenant compte de leurs multiplicit\'es, les valeurs propres de $D_M$ par $\lambda_k$ ($k\geq 1$), et supposerons la suite $(|\lambda_k|)_{k\geq 1}$ croissante.\\ Af\mbox{}in d'estimer ce spectre en fonction d'invariants extrins\`eques, nous comparons $D$ \`a un op\'erateur ici auxiliaire appel\'e \emph{op\'erateur de Dirac-Witten} \cite{Wi81,Zhang98} et d\'ef\mbox{}ini dans une base orthonorm\'ee locale $(e_j)_{1\leq j\leq m}$ de $TM$ par $\wih{D}:=\sum_{j=1}^me_j\cdot\wnabla_{e_j}$. Gr\^ace aux propri\'et\'es de l'isomorphisme (\ref{isomspin}), les op\'erateurs $D^2$ et $\wih{D}^2$ sont li\'es par (cf. \cite{Ginthese} ou \cite{Gin2002}, Lemma 2):
\begin{equation}\label{lcarresDwihD}
\wih{D}^2\phi=D^2\phi-\frac{m}{2}dH\hspace{-1mm}\cdotM\hspace{-0.5mm}\phi-\frac{m^2H^2}{4}\phi,\end{equation} identit\'e valable pour toute section $\phi$ de $\Sigma\wit{M}_{|_M}$.\\

\section{R\'esultat principal}
\noindent Nous nous restreignons maintenant au cas o\`u la vari\'et\'e ambiante $(\wit{M},g)$ est l'espace hyperbolique r\'eel $\mathbb{H}^{m+1}$ muni de sa m\'etrique standard $g$ \`a courbure sectionnelle $-1$. En tant que vari\'et\'e riemannienne spinorielle, l'espace hyperbolique poss\`ede la propri\'et\'e remarquable d'admettre des \emph{spineurs de Killing imaginaires}, i.e., il existe des sections non nulles $\phi$ de $\Sigma\wit{M}$ satisfaisant, pour tout champ de vecteurs $X$,
\begin{equation}\label{eqspKI}\wnabla_X\phi=\pm\frac{i}{2}X\cdot\phi.\end{equation} 
Il est ici \`a noter que la classif\mbox{}ication des vari\'et\'es riemanniennes spinorielles compl\`etes admettant de telles sections a \'et\'e achev\'ee par H. Baum dans \cite{Baum892} et \cite{Baum891}, faisant appara\^\i tre d'autres exemples que l'espace hyperbolique. Pour les propri\'et\'es des spineurs de Killing, on pourra consulter \cite{BFGK}.\\
Reprenant l'id\'ee introduite dans \cite{Baer98} d'utiliser la restriction de ces spineurs de Killing comme sections-test dans le principe du Min-Max, nous poussons plus loin la preuve du r\'esultat (3) de C. B\"ar gr\^ace \`a l'identit\'e (\ref{lcarresDwihD}) qui permet d'\'eviter l'emploi de l'in\'egalit\'e de Cauchy-Schwarz (comparer avec \cite{Baer98}, p. 586). Nous d\'emontrons le th\'eor\`eme suivant:\\
\begin{thm}\label{tmajKIHyp}
Soit $(M,\,g)$ une hypersurface riemannienne (immerg\'ee) de dimension $m$
compacte et orient\'ee de l'espace hyperbolique $(\mathbb{H}^{m+1},g)$. Munissons $M$ de la structure spinorielle induite et posons $N:=2^{[\frac{m+2}{2}]}$. 
Alors pour tout $k\in\{1,\ldots,N\}$,
\begin{equation}\label{majKIHyp}
\lambda_k^2\leq\frac{m^2}{4}\pa{\buil{\sup}{M}\hspace{-0.5mm}H^2-1}.
\end{equation} De plus, si (\ref{majKIHyp}) est une \'egalit\'e pour $k=1$, alors la courbure moyenne $H$ est constante.\\
\end{thm}

\noindent{\it D\'emonstration}: \'Etant donn\'ee une section non nulle $\phi$ de $\Sigma\mathbb{H}^{m+1}$ satisfaisant (\ref{eqspKI}), \'evaluons le quotient de Rayleigh
$\mathcal{Q}\pa{D^2,\phi}:=\frac{\int_M<D^2\phi,\phi>v_g}{\int_M<\phi,\phi>v_g}$. La section $\phi$ v\'erif\mbox{}iant (\ref{eqspKI}), il vient imm\'ediatement $\wih{D}\phi=\mp \frac{mi}{2}\phi$ puis $\wih{D}^2\phi=-\frac{m^2}{4}\phi$. D'apr\`es (\ref{lcarresDwihD}), 
\[D^2\phi=-\frac{m^2}{4}\phi+\frac{m^2H^2}{4}\phi+\frac{m}{2}dH\hspace{-1mm}\cdotM\hspace{-0.5mm}\phi.\] Prenons le produit scalaire hermitien de cette identit\'e avec $\phi$ et identif\mbox{}ions-en les parties r\'eelles: puisque la multiplication de Clif\mbox{}ford par une $1$-forme est anti-autoadjointe (cf. \cite{LM}), $\Re\textrm{e}(<dH\hspace{-1mm}\cdotM\hspace{-0.5mm}\phi,\phi>)=0$, dont on d\'eduit que:
\begin{equation}\label{D2KI}\Re\textrm{e}\pa{<D^2\phi,\phi>}=-\frac{m^2}{4}<\phi,\phi>+\frac{m^2H^2}{4}<\phi,\phi>.\end{equation} Par int\'egration de cette identit\'e sur $M$, nous obtenons
\[
\mathcal{Q}\pa{D^2,\phi}=-\frac{m^2}{4}+\frac{m^2\int_MH^2<\phi,\phi>v_g}{4\int_M <\phi,\phi>v_g}\leq-\frac{m^2}{4}+\frac{m^2}{4}\buil{\sup}{M}\hspace{-0.5mm}H^2.\] Puisque l'espace hyperbolique admet un espace de dimension $2.2^{[\frac{m+1}{2}]}$ de spineurs de Killing imaginaires, l'application du principe du Min-Max donne l'in\'egalit\'e (\ref{majKIHyp}).\\
Si maintenant (\ref{majKIHyp}) est une \'egalit\'e pour $k=1$, alors la caract\'erisation variationnelle de $\lambda_1^2$ entraine, pour la restriction \`a $M$ de tout spineur de Killing imaginaire $\phi$ sur $\mathbb{H}^{m+1}$,
\[D^2\phi=\lambda_1^2\phi.\] L'injection de cette relation dans (\ref{D2KI}) et le fait qu'un spineur de Killing imaginaire n'admet pas de z\'ero \cite{BFGK} donne alors $\lambda_1^2=\frac{m^2}{4}\pa{H^2-1}$, dont on d\'eduit que $H$ doit \^etre constante.\findemo
{\bf Remarques}\\
\noindent a. Ce r\'esultat se r\'ev\`ele analogue \`a celui de E. Heintze (\cite{Hein88}, Theorem 2.3) pour le laplacien scalaire, et est optimal: pour toute sph\`ere g\'eod\'esique, (\ref{majKIHyp}) est une \'egalit\'e pour $k=1$ (voir \cite{Gin2002}).\\
b. L'in\'egalit\'e (\ref{majKIHyp}) pour $k=1$ demeure \'evidemment valable sur toute hypersurface (compacte et orient\'ee) de toute vari\'et\'e riemannienne spinorielle admettant un spineur de Killing imaginaire.\\
c. Pour obtenir un majorant $L^2$ analogue \`a ceux de (1) et (2), il suf\mbox{}f\mbox{}irait qu'il existe un spineur de Killing imaginaire de norme \emph{constante} sur $M$. Il a n\'eanmoins \'et\'e d\'emontr\'e \cite{Baum892,Baum891,Ginthese} que, sous cette hypoth\`ese, $M$ doit \^etre totalement ombilique dans $\mathbb{H}^{m+1}$, i.e., \^etre isom\'etrique \`a une sph\`ere g\'eod\'esique de $\mathbb{H}^{m+1}$.\\

\noindent Nous n'excluons pas l'existence lorsque $\wit{M}=\mathbb{H}^{m+1}$ d'une telle borne $L^2$, qui para\^\i trait naturelle au regard de l'estimation obtenue par A. El Souf\mbox{}i et S. Ilias (\cite{ElSI92}, Th\'eor\`eme 1) pour la premi\`ere valeur propre non nulle du laplacien scalaire. Cependant, nous signalons les limites de cette analogie: l'estimation \emph{a priori} (comparer avec \cite{ElSI92}, Th\'eor\`eme 2)
\[\lambda_1^2\leq\frac{m^2}{4\mathrm{Vol}(M)}\int_M\pa{H^2+R(\iota)}v_g\] sur une hypersurface compacte d'une quelconque vari\'et\'e $\wit{M}$ \emph{conform\'ement} immerg\'ee dans la sph\`ere $S^{m+1}$ ne peut avoir lieu en dimension $m=2$: en ef\mbox{}fet, la quantit\'e $\int_M\pa{H^2+R(\iota)}v_g$ est dans ce cas la fonctionnelle de Willmore, qui est invariante par changement conforme de m\'etrique sur $\wit{M}$. Or le produit $\lambda_1^2\mathrm{Aire}(S^2)$ est \emph{non born\'e} dans la classe conforme de $\mathrm{can}_{S^2}$ \cite{Amm02}. Cette dif\mbox{}f\'erence de comportement entre l'op\'erateur de Dirac et le laplacien scalaire se retrouve d'ailleurs dans la caract\'erisation du cas d'\'egalit\'e de (2) \cite{GinEgBaer2003}.\\ Pour l'utilisation des propri\'et\'es conformes de l'op\'erateur de Dirac dans ce contexte, nous renvoyons \`a \cite{Gin2002,Ginthese}.

\providecommand{\bysame}{\leavevmode\hbox to3em{\hrulefill}\thinspace}

$ $\\
\noindent \begin{small}Nicolas GINOUX\\
Max-Planck Institut f\"ur Mathematik in den Naturwissenschaften, Inselstra\ss e 22 D-04103 Leipzig.\\ ginoux@mis.mpg.de \end{small}
\end{document}